\documentclass[largeformat,article,a4paper]{interact}

\usepackage{graphicx} 
\usepackage{multicol}        % used for the two-column index
\usepackage{amsmath,amssymb}
\usepackage{mathtools}
\usepackage{bm}
\usepackage{esvect}
\usepackage{mathptmx}      % use Times fonts if available on your TeX system
\usepackage{latexsym}
\usepackage{cancel}
\usepackage{color}
\usepackage{multirow}
\usepackage[numbers]{natbib}
\usepackage{tikz}
\usepackage{tikzscale}
\usepackage{pgfplots}
\usepackage{soul} %per barrare testo con \st
\graphicspath{{./Figures/}} 

\usetikzlibrary{shapes,arrows}
\usetikzlibrary{calc,matrix,decorations.markings,decorations.pathreplacing}

\newcommand{\translatepoint}[1]%
{   \coordinate (mytranslation) at (#1);
}

\tikzset{ 
  table/.style={
    matrix of nodes,
    row sep=-\pgflinewidth,
    column sep=-\pgflinewidth,
    nodes={rectangle,text width=2cm,align=center},
    text depth=1.25ex,
    text height=2.5ex,
    nodes in empty cells
  }
}

\newcommand{\bu}{\mathbf{u}}
\newcommand{\bv}{\mathbf{v}}
\newcommand{\bC}{\mathbf{C}}

\newcommand{\bV}{\mathbf{V}}
\newcommand{\bO}{\mathbf{S}}
\newcommand{\Cop}{\boldsymbol{\mathcal{C}}}
\newcommand{\Lop}{\boldsymbol{\mathfrak{L}}}
\newcommand{\bD}{\mathbf{D}}
\newcommand{\bU}{\mathbf{U}}

\newcommand{\bM}{\mathbf{M}}

\newcommand{\bG}{\mathbf{G}}
\newcommand{\bcA}{\boldsymbol{\mathcal{A}}}
\newcommand{\bcD}{\boldsymbol{\mathcal{D}}}

\newcommand{\bui}{\mathbf{u}_i}
\newcommand{\bun}{\mathbf{u}^n}
\newcommand{\bunp}{\mathbf{u}^{n+1}}

\newcommand{\buj}{\mathbf{u}_j}

\newcommand{\Ec}{K_e}
\newcommand{\Hc}{K_h}

\begin{document}

%\articletype{ARTICLE TEMPLATE}
\articletype{}

\title{Effects of discrete energy and helicity conservation in numerical simulations of helical turbulence}

\author{
\name{Francesco Capuano\textsuperscript{a} and Donato Vallefuoco\textsuperscript{b}}
 \affil{\textsuperscript{a}Dipartimento di Ingegneria Industriale (DII), Universit\`a di Napoli Federico II, Napoli, 80125, Italy}
\affil{\textsuperscript{b}Laboratoire de M\'ecanique des Fluides et d'Acoustique, CNRS, \'Ecole Centrale de Lyon, Universit\'e Claude Bernard Lyon 1 and INSA de Lyon, 36 avenue Guy de Collongue, F-69134 \'Ecully Cedex, France}
}

\maketitle

\begin{abstract}
Helicity is the scalar product between velocity and vorticity and, just like energy, its integral is an inviscid invariant of the three-dimensional incompressible Navier-Stokes equations. However, space- and time-discretization methods typically corrupt this property, leading to violation of the inviscid conservation principles. This work investigates the discrete helicity conservation properties of spectral and finite-differencing methods, in relation to the form employed for the convective term. Effects due to Runge-Kutta time-advancement schemes are also taken into consideration in the analysis. 
The theoretical results are proved against inviscid numerical simulations, while
a scale-dependent analysis of energy, helicity and their non-linear transfers is performed to further characterize the discretization errors of the different forms in forced helical turbulence simulations. %oppure: 
%decaying helical turbulence simulations are performed to further characterize the effect of helicity conservation.
\end{abstract}

\section{Introduction} \label{intro}
This work is concerned with the incompressible Navier-Stokes (NS) equations,
which can be written as follows
%
%\begin{eqnarray}
%\frac{\partial {u}_{i}}{\partial x_{i}}=0 \;, \label{eq:NavSto1} \\
%\frac{\partial {u}_{i}}{\partial t} + \mathcal{N}_{i}({u}) = -\frac{\partial {p}}{\partial x_{i}} + \frac{1}{\textrm{Re}}\frac{\partial^2 {u}_{i}}{\partial x_{j}\partial x_{j}}
%\;, \label{eq:NavSto2}
%\end{eqnarray}
%
%
\begin{eqnarray}
\nabla \cdot \vv{u}= 0 \;, \label{eq:NavSto1} \\
\frac{\partial \vv{u}}{\partial t} + \vv{\mathcal{N}}(\vv{u}) = -\nabla p + \nu\nabla^2 \vv{u}
\;, \label{eq:NavSto2}
\end{eqnarray}
where $\vv{\mathcal{N}}(\vv{u})$ is the non-linear convective term and $\nu$ is the kinematic viscosity.
%
%In the inviscid limit of $\text{Re} \rightarrow \infty$, the NS equations
In the inviscid case the NS equations are replaced by the Euler equations, which
possess \textit{linear} as well as \textit{quadratic} invariants, i.e., linear or quadratic functions of the velocity vector whose integral over the domain remains constant in time throughout the flow evolution. The identification of the invariants becomes particularly straightforward upon recognizing that the inviscid NS equations are a Hamiltonian system and can be thus written in the canonical skew-gradient formulation
\begin{equation}
\dfrac{\partial \vv{y}}{\partial t} = \mathcal{D} \left(\vv{y} \right) \nabla \mathcal{H} \;, \label{eq:skewgrad}
\end{equation}
where $\mathcal{D}$ is a skew-symmetric operator and $\mathcal{H}$ is the Hamiltonian functional \cite{olver1982nonlinear}. Recasting the inviscid version of Eqs.~(\ref{eq:NavSto1})-(\ref{eq:NavSto2}) in the skew-gradient formulation of Eq.~(\ref{eq:skewgrad}) involves either eliminating the pressure (upon projection of the velocity field onto a divergence-free space) \cite{arnold1969hamiltonian}, or taking the curl of Eq.~(\ref{eq:NavSto2}) \cite{olver1982nonlinear}.
Once a system is known to be Hamiltonian, the Noether theorem can be applied to yield an association between the symmetries of the Euler system and conserved quantities; for instance, conservation of linear momentum corresponds to the invariance of the Euler equations under space-time translations \cite{olver1980hamiltonian}. Here we are particularly interested in quadratic invariants. One follows directly from Eq.~(\ref{eq:skewgrad}), since the Hamiltonian functional is seen to be the mean kinetic energy of the flow, 
\begin{equation}
\Ec = \dfrac{1}{2} \langle \vv{u} \cdot \vv{u} \rangle \;,
\end{equation}
where the angular brackets denote the average operator over the entire domain.
Another inviscid integral arises from the nontrivial kernel of the operator $\mathcal{D}$,
and corresponds to the quantity named \textit{helicity}, whose average value over the domain reads
\begin{equation}
\Hc = \langle \vv{u} \cdot \vv{\omega} \rangle \;,
\end{equation}
where $\vv{\omega} = \nabla \times \vv{u}$ is the vorticity vector. Helicity has thus been termed a \textit{Casimir} invariant of the Euler equations \cite{holm2011geometric}.

The discovery of helicity conservation is relatively recent and dates back to 1961 \cite{moreau1961constantes}. Since then, helicity has been found to play an important role in laminar and turbulent flows, and has raised significant interest both in terms of fundamental understading and repercussions in engineering applications and flows of practical interest \cite{moffatt1992helicity}.
Helicity is known to inhibit the transfer of energy towards smaller scales, since the statistical alignment of velocity and vorticity leads to partial suppression of the nonlinear term \cite{andre1977influence,kraichnan1988depression}. On the other hand, the self-similar energy decay rate has been shown to be the same as in the nonhelical case, unless rotation is present \cite{teitelbaum2009effect}. The dynamics of the joint cascade of energy and helicity constitutes an interesting field of research which is very active nowadays \cite{chen2003joint,biferale2013split,kessar2015non}.
Helicity is also important in atmospheric and geophysical flows \cite{lautenschlager1988subgrid}, while an helicity-based index is used in biomedical research to quantify swirling motions in cardiovascular flows \cite{morbiducci2007helical}.

Upon discretization of the NS equations in space and time, the invariant character of energy and helicity is generally lost. A significant research effort has been carried out over the last years to develop numerical algorithms that preserve invariants also in a discrete sense, with the aim of obtaining stable computations and physically relevant solutions. 
Particular attention has been paid to the development of energy-preserving numerical methods, that have ultimately allowed stable long-time integrations and realistic representations of the energy cascade \cite{Morinishi1998,Verstappen2003,Capuano2015b}. On the other side, invariance of helicity has been very seldom considered in the derivation of numerical methods. Notable exceptions include the works by Liu and Wang \cite{liu2004energy}, for axisymmetric flows, and by Rebholz and coworkers (see, e.g., \cite{rebholz2007energy}), in the framework of finite-element methods.

The aim of the present work is to investigate the discrete helicity-conservation properties 
of commonly used semi-discrete algorithms based on spectral, finite-difference or finite-volume methods and Runge-Kutta time-advancement schemes. Particularly, the discrete behaviour of the different formulations of the nonlinear term $\vv{\mathcal{N}}$ is discussed. The effects of discretization in time are also taken into consideration in the analysis.

The paper is organized as follows. Section~\ref{sec:1} deals with the energy and helicity balance equations in a continuous setting, along with the rules of calculus required to derive them. The spatial and temporal discretization of the Navier-Stokes equations are presented in Section~\ref{sec:2}, in which a fully discrete equation for the helicity evolution is also derived. 
In Section~\ref{sec:results} the numerical results are reported and discussed. Concluding remarks are given in Section~\ref{sec:conclusions}.

\section{Conservation of energy and helicity} \label{sec:1}

In this section, we derive the conservation laws for the two quadratic invariants of the Euler equations. Clearly, in the viscous case, these quantities are actually not conserved but are subject to a source term which is either strictly disspative (for the energy) or has no definite sign (for the helicity).

For the sake of clarity, we will assume that periodic boundary conditions apply, but this does not come at a loss of generality. In deriving the conservation laws, it will
be useful to list a number of identities, which can be readily proven using the 
standard rules of calculus 
\begin{eqnarray}
\nabla \cdot (\vv{u}\vv{u}) &=& \vv{u} \cdot \nabla \vv{u} \; \quad \text{if} \; \nabla \cdot \vv{u} = 0 \;, \label{eq:prodrule} \\
\langle \nabla u \rangle &=& \vv{0} \;, \label{eq:divrule} \\
\langle \left(\nabla u\right) v \rangle &=& - \langle u \left(\nabla v\right) \rangle \;, \label{eq:sumbyparts} \\
\langle \vv{u} \cdot \left(\nabla \times \vv{v} \right) \rangle &=& \langle \vv{v} \cdot \left(\nabla \times \vv{u} \right) \rangle \label{eq:curl} \;.
%\langle \left(\partial_i f\right) g \rangle
\end{eqnarray}
All the functions appearing above are supposed to be periodic.

The energy evolution equation can be derived by scalarly multiplying Eq.~(\ref{eq:NavSto2}) by $\vv{u}$ and taking the average over the domain, yielding
\begin{equation}
\dfrac{\text{d} \Ec}{\text{d} t} + \langle \vv{u} \cdot \vv{\mathcal{N}} \rangle = 
- \cancel{\langle \vv{u} \cdot \nabla p \rangle} + \nu \langle \vv{u} \cdot \nabla^2 \vv{u} \rangle \;. \label{eq:energy}
\end{equation}
The contribution from the pressure term vanishes due to use of the integration-by-parts~(\ref{eq:sumbyparts}) together with the incompressibility condition~(\ref{eq:NavSto1}).
Similarly, for helicity one also has to take the curl of the momentum equation~(\ref{eq:NavSto2}) and then scalarly multiply the resulting vorticity equation by $\vv{u}$, leading to
\begin{equation}
\dfrac{1}{2} \dfrac{\text{d} \Hc}{\text{d} t} + \langle \vv{u} \cdot \nabla \times \vv{\mathcal{N}} \rangle =  \nu \langle \vv{u} \cdot \nabla^2 \vv{\omega} \rangle \;. \label{eq:helicity}
\end{equation}
The contributions from the nonlinear term in both Eq.~(\ref{eq:energy}) and Eq.~(\ref{eq:helicity}) also vanish, but the properties to be invoked to prove their cancellation depend on the formulation adopted for the nonlinear term. This can be indeed expressed in several, analytically equivalent forms:
\begin{eqnarray}
(\text{Adv.})&  \quad \vv{\mathcal{N}} &= \vv{u} \cdot {\nabla \vv{u}} \;, \\
(\text{Div.})&  \quad \vv{\mathcal{N}} &= \nabla \cdot (\vv{u} \vv{u}) \;, \\
(\text{Skew.})& \quad \vv{\mathcal{N}} &= \dfrac{1}{2} \left(\vv{u} \cdot {\nabla \vv{u}} + \nabla \cdot (\vv{u} \vv{u}) \right) \;, \\
(\text{Rot.})&  \quad \vv{\mathcal{N}} &= \vv{\omega} \times \vv{u} \;,
\end{eqnarray}
which are named the advective, divergence, skew-symmetric and rotational forms respectively \cite{Capuano2015b}. Switching from one form to another requires use of the product rule \eqref{eq:prodrule} and of Eq.~(\ref{eq:NavSto1}). Note that when the rotational form is used, the kinematic pressure is subsituted by the dynamic pressure $P = p + |\vv{u}|^2/2$. The properties required to prove cancellation of the nonlinear term contribution in the energy and helicity equations are reported in Table~\ref{table:properties}.
Note that while the advective and divergence forms always require the use of the product rule, only the integration-by-parts property is needed to prove energy conservation for the skew-symmetric form. However, the skew-symmetric form still necessitates the product rule to verify helicity invariance.
On the other hand, the rotational form is automatically energy conserving (due to orthogonality between the velocity and its cross product with vorticity), and requires only the symmetry of the curl operator to demonstrate conservation of helicity.
%It is worth to note that the integration-by-parts rule is still required to prove cancellation of the pressure term in the energy equation} \textcolor{red}{credo non sia attinente, e un pezzo che non fa piu parte del termine non lineare. Inoltre, nel caso incompressibile, la pressione e' un moltiplicatore di Lagrange, se ci aggiungiamo il pezzo a gradiente non cambia nulla ne' nel continuo ne' nel discreto}.

The energy and helicity equations can be then integrated in time between two generic instants $t$ and $t+\Delta t$, to yield the finite time-variation of the two invariants
\begin{eqnarray}
\dfrac{\Delta \Ec}{\Delta t} = \left\langle \nu \langle \vv{u} \cdot \nabla^2 \vv{u} \rangle \right\rangle_{\Delta t} \;, \label{eq:energy_int} \\[0.2cm]
\dfrac{\Delta \Hc}{\Delta t} = \left\langle 2\nu \langle \vv{u} \cdot \nabla^2 \vv{\omega} \rangle \right\rangle_{\Delta t} \label{eq:helicity_int} \;,
\end{eqnarray}
where $\langle \cdot \rangle_{\Delta t}$ is the time-average operator over the interval $\Delta t$.
In deriving Eqs.~(\ref{eq:energy_int})-(\ref{eq:helicity_int}), which will be useful in the following, we have implicitly employed the integration-by-parts rule for the time-derivative operator. Indeed, given two differential equations $\partial_t f(t) = F$ and $\partial_t g(t) = G$, where $f(t)$ and $g(t)$ are two continuous and differentiable functions of time, and $\partial_t$ is a continuous time-derivative operator, the integration-by-parts rule gives
\begin{equation}
	\dfrac{\Delta (f \cdot g)}{\Delta t} = \left\langle  f \cdot G \right\rangle_{\Delta t}  + \left\langle g \cdot F \right\rangle_{\Delta t} \;, \label{eq:SBPtime_cont}
\end{equation}
which is needed to prove Eqs.~(\ref{eq:energy_int})-(\ref{eq:helicity_int}) upon substituting $f$ and $g$ with $\vv{u}$, and with $\vv{u}$ and $\vv{\omega}$ respectively.

The properties (\ref{eq:prodrule})-(\ref{eq:curl}) and (\ref{eq:SBPtime_cont}) do not necessarily hold on a discrete level when the continuous space and time operators are substituted by their finite-differencing counterparts, as will be outlined in the following section.

\begin{table}
\centering
  \caption{Minimum set of properties required to prove cancellation of the nonlinear term contribution in the energy and helicity equations for the various formulations of convection in a continuous setting. \label{table:properties}}
  \begin{tabular}{llllllllll}
   \hline\noalign{\smallskip}
    & &  \multicolumn{4}{l}{Energy} & \multicolumn{4}{l}{Helicity} \\
   %\cmidrule{2-9}
   Eq. & Property & Adv. & Div. & Skew. & Rot. & Adv. & Div. & Skew. & Rot. \\
   \noalign{\smallskip}\hline\noalign{\smallskip}
   (\ref{eq:prodrule})   & Product rule            & $\times$ & $\times$ &   &  & $\times$ & $\times$ & $\times$ &   \\
   (\ref{eq:divrule})    & Divergence              & $\times$ & $\times$ &   &  &  &  &  &   \\
   (\ref{eq:sumbyparts}) & Integration by parts    &   &   & $\times$ &  &  &  &  &   \\
   (\ref{eq:curl})       & Symmetric curl          &   &   &   &  & $\times$ & $\times$ & $\times$ & $\times$ \\
   \noalign{\smallskip}\hline
  \end{tabular}
\end{table}

\section{Fully discrete energy and helicity evolution} \label{sec:2}

The ultimate aim of this section is to derive a fully discrete analogue of Eqs.~(\ref{eq:energy_int})-(\ref{eq:helicity_int}). In doing this, we also wish to understand which of the properties listed in Table~\ref{table:properties}, as well as Eq.~(\ref{eq:SBPtime_cont}) are retained by the discretization procedure.
To this end, we employ a semi-discrete approach (also called \textit{method of lines}) which is routinely used to treat the incompressible Navier-Stokes equations numerically. 

\subsection{Spatial discretization} \label{sec:spatial}

The first step of the algorithm is to derive a spatially discretized version of Eqs.~(\ref{eq:NavSto1})-(\ref{eq:NavSto2}), which can be expressed as follows
\begin{eqnarray}
\mathbf{M}\bu = \mathbf{0} \;, \label{eq:NavStoDiscr1} \\
\frac{\text{d} \bu}{\text{d} t} + \bC(\bu)\bu = -\mathbf{G}\mathbf{p}
+ \nu \mathbf{L}\bu\;, \label{eq:NavStoDiscr2}
\end{eqnarray}
where $\bu$ is a discrete vector containing the
components of velocity on the three-dimensional mesh,
$\bu = \left[ \bu_x \; \bu_y \; \bu_z \right]^T$,
the matrices $\mathbf{G}\in {R}^{N_{\bu}\times N_p}$ and
$\mathbf{M}\in {R}^{N_p\times N_{\bu}}$ are the discrete
gradient and divergence operators, respectively, while
$\mathbf{L}\in {R}^{N_{\bu}\times N_{\bu}}$ is the
block-diagonal Laplacian, $\text{diag}(\Lop,\Lop,\Lop)$.
%with $\Lop\in R^{N_p\times N_{p}}$.
$N_\bu$ and $N_p$ are the number of unknowns on the
mesh for velocity and pressure respectively.
The spatial discretization reported in Eqs.~(\ref{eq:NavStoDiscr1})-(\ref{eq:NavStoDiscr2})
possibly encompasses finite-difference, finite-volume, pseudo-spectral
as well as finite-element approaches, although the details of the 
various operators depend on the specific method employed.
In what follows, it is assumed that the
discretization is built upon a regular mesh
with uniform step size, and that gradient and divergence
operators are discretized consistently, in such a way that the
relation $\bG^T=-\bM$ holds. The uniform mesh assumption
is made here only for the sake of clarity; extension to more general cases
is discussed in Section~\ref{sec:irregular}.

We are particularly interested in the convective term, 
which can be expressed as the product of a linear
block-diagonal convective operator $\bC(\bu)$ and $\bu$:
\begin{equation}
\bC(\bu) \bu = \left[\begin{array}{ccc}
 \Cop_x(\bu) & & \\
 & \Cop_y(\bu) & \\
 & & \Cop_z(\bu) \\
 \end{array}\right]
\left[\begin{array}{c}
 \bu_x \\
 \bu_y \\
 \bu_z \\
 \end{array}\right]
  \; .
\end{equation}
The operator $\bC(\bu)$ is obtained by discretizing the nonlinear term
starting from one of the possible expressions in which this can be written.
As mentioned in Section~\ref{sec:1}, while all these expressions
are equivalent in a continuous framework, this is not the case
for the discrete setting.
Also, the operator $\Cop$ assumes a specific form which depends on the
layout of the variables on the mesh. Typically, 
the Navier-Stokes equations are either discretized in a 
\textit{colocated} arrangement, i.e., with all the variables
located at the same points of the mesh, or on a \textit{staggered} grid,
in which the velocity components and the pressures are stored
in different locations.

In the subsequent sections we will first 	analyze the colocated layout of the variables,
and then discuss the staggered one.

\subsubsection{Colocated layout}

In a colocated layout, the operators $\Cop_x$, $\Cop_y$ and $\Cop_z$ are identical
and can assume one of the following forms:
\begin{eqnarray}{\label{eq:NSCOP} }
(\text{Adv.})&  \quad \Cop(\bu)  &= \left( \bU_x \bD_x + \bU_y \bD_y +
   \bU_z \bD_z \right) \equiv \bcA \;, \label{eq:NSCOP_adv} \\
(\text{Div.})&  \quad \Cop(\bu)  &= \left( \bD_x \bU_x + \bD_y \bU_y +
   \bD_z \bU_z \right) \equiv \bcD \;, \label{eq:NSCOP_div} \\
(\text{Skew.})& \quad \Cop(\bu) &=  \frac{1}{2}\left(\bcA + \bcD \right) \;,
\label{eq:NSCOP_skew}
\end{eqnarray}
while the rotational form can be expressed as
\begin{equation}
(\text{Rot.}) \quad \bC(\bu) = \mathbf{V}\left(\mathbf{R} \bu \right) = - \mathbf{V}\left( \bu \right)\mathbf{R}  \;. \label{eq:NSCOP_rot}
\end{equation}
In the expressions above, the matrices $\bD_{(\cdot)}$
represent the discrete derivative
operators along each direction and acting on the whole
set of variables on the mesh, while
$\bU_{(\cdot)}$ are the diagonal matrices of the discretized
velocity components along the three directions
(e.g., $\bU_x = \text{diag}(\bu_x)$).
Also, $\mathbf{V}$ is a skew-symmetric matrix performing pointwise vector product, 
while $\mathbf{R}$ is the curl operator, which can be expressed as
\begin{equation}
\mathbf{R} = \left[\begin{array}{ccc}
 & -\mathbf{D}_z & \mathbf{D}_y \\
 \mathbf{D}_z &  & -\mathbf{D}_x \\
 -\mathbf{D}_y & \mathbf{D}_x & \\
 \end{array}\right] \;. \label{eq:curl_discrete}
\end{equation}

Generally speaking, the introduced discretization yields truncation errors, due to the finite approximations
of derivatives, as well as aliasing-type errors, due to the evaluation of the nonlinear terms in
a finite-dimensional space.
These errors are likely to produce deviations from the rules of calculus and therefore
violation of the conservation of quadratic invariants. 
%\textcolor{red}{However, in dealiased pseudospectral simulations these errors vanish and quadratic invariants are conserved, the only numerical approximation being the spectral truncation which does not alter energy and helicity conservation.}
As a first step, it is useful to check how the discrete operators behave
with respect to the properties listed in Table~\ref{table:properties}.
To this end, we introduce a discrete inner norm, which under the hypothesis
of regular mesh simply reads
\begin{equation}
\langle u  v \rangle \quad \rightarrow \quad \bu^T \bv.
\end{equation}
The discrete integration-by-parts rule is easily verified for central-differencing
formulas, including compact schemes and Fourier differentiation, which all 
lead to skew-symmetric derivative matrices. Indeed, in such cases, assuming two discrete 
functions $\bu$ and $\bv$ on a periodic domain, and a generic derivative operator $\mathbf{D}$, one has 
\begin{equation} 
\bu^T\bD\bv = \bv^T\bD^T\bu = -\bv^T\bD\bu \;, \label{eq:SBP_discr}
\end{equation}
which is the discrete analogue of Eq.~(\ref{eq:sumbyparts}).
It is worth to note that Eq.~(\ref{eq:SBP_discr}) holds regardless
of the scheme order (i.e., the truncation error), as well as of the presence of aliasing errors.
%which are indeed present in the pointwise evaluation of the products in Eq.~(\ref{eq:SBP_discr}).
%\textcolor{red}{per lo stesso discorso di prima credo che l'aliasing dipenda dalla proiezione sulla base di Fourier. che ne dici di '...regardless
%of the scheme order (i.e., the truncation error), as well as of aliasing errors (in spectral aliased simulations).'}
This was already shown by Kravchenko and Moin \cite{Kravchenko1997} using Fourier expansion.
The divergence property is also easily proven assuming that the discrete derivative 
operators behave consistently with respect to a constant function. In this case, 
one has simply
\begin{equation} 
\mathbf{1}^T\bD\bu = -\bu^T\bD\mathbf{1} = 0 \;, \label{eq:div_discr}
\end{equation}
where $\mathbf{1}$ is the unit vector.
As a consequence of using central-differencing schemes, the curl operator
reported in Eq.~(\ref{eq:curl_discrete})
also turns out to be symmetric as in the continuous case, Eq.~(\ref{eq:curl}),
being a skew-symmetric block matrix constituted by skew-symmetric blocks.
On the other hand, it is straightforward to verify that the product rule is in general violated by 
the discrete approximation. %\textcolor{red}{questa frase dice lo stesso della frase successiva, oppure non ho capito quest'ultima frase a cosa si riferisce}
Although this property is restored in case of dealiased spectral
schemes, for finite-differencing methods of any order the product rule does not
hold on a discrete level, even if aliasing errors are removed \cite{Kravchenko1997}.
A summary of the properties is reported in Table~\ref{table:properties_discr}.

\begin{table}
\centering
  \caption{Requirements for a spatial discretization on a colocated grid to satisfy the rules of calculus. \label{table:properties_discr}}
  \begin{tabular}{lll}
   \hline\noalign{\smallskip}
   Eq. & Property & Requirement \\
   \noalign{\smallskip}\hline\noalign{\smallskip}
   (\ref{eq:prodrule})   & Product rule            & Fully de-aliased spectral differentiation \\
   (\ref{eq:divrule})    & Divergence              & Consistent scheme ($\mathbf{D} \mathbf{1} = 0$) \\
   (\ref{eq:sumbyparts}) & Integration by parts    & \multirow{ 2}{*}{Central-differencing scheme (incl. explicit, compact, spectral)} \\
   (\ref{eq:curl})       & Symmetric curl          &  \\
   \noalign{\smallskip}\hline
  \end{tabular}
\end{table}

In light of the above considerations, and upon comparison with Table~\ref{table:properties},
one could already infer the consequences of using the various formulations
of the nonlinear term. Further insight, however, is gained by deriving
the semi-discrete, time-continuous energy and helicity equations.
The discrete energy is defined as $e = \bu^T \bu/2$, whereas the discrete helicity
reads $h = \mathbf{u}^T \mathbf{w}$, where $\mathbf{w} = \mathbf{R} \bu$
is the discrete vorticity vector. The energy equation has been widely
reported elsewhere \cite{Capuano2015b,Verstappen2003} and reads
\begin{equation}
\dfrac{\text{d} e}{\text{d} t}  =  -\bu^T\mathbf{C}(\bu) \bu -
\bu^T \mathbf{G}\mathbf{p} + \nu \bu^T \mathbf{L} \bu \;. \label{eq:energy_semidiscr}
\end{equation}
While the pressure term vanishes as a consequence of the continuity equation
(assuming $\mathbf{M} \bu = \mathbf{0}$ and the consistency relation $\bG^T=-\bM$), the convective term contribution
depends basically on the formulation employed for it. 
It is easy to show that Eq.~(\ref{eq:NSCOP_skew}) yields a skew-symmetric 
convective operator provided that the schemes satisfy the discrete integration-by-parts
rule, while the expression~(\ref{eq:NSCOP_rot}) is skew-symmetric a priori.
On the other hand, this does not happen for none of the operators 
given by Eq.~(\ref{eq:NSCOP_adv}) and Eq.~(\ref{eq:NSCOP_div}). 
As a consequence, only the skew-symmetric and the rotational
forms of the convective term preserve energy spatially in the inviscid case.
This circumstance was soon recognized \cite{Zang1991}, and since then 
there has been a long debate regarding the use of these two forms,
which allow stable long-time integration \cite{Blaisdell1996,Horiuti}. 
While the rotational form is much more cost-effective
than the skew-symmetric one, it is generally agreed that
its aliasing errors are higher than those provided
by the skew-symmetric form \cite{Zang1991}. Indeed, aliased simulations carried out with the skew-symmetric form
proved to be very similar to the de-aliased ones \cite{Kravchenko1997}.
%The skew-symmetric form is nowadays regularly applied in several highly accurate finite-difference and spectral numerical codes for DNS and LES of turbulent flows 
%(\cite{Laizet2009}-\cite{Gibson2014}). The conservation properties
%of the rotational form are also massively exploited 
%(\cite{Mahesh2004}-\cite{Vallefuoco2017}).

Here, we are especially interested in how the helicity balance
is affected by the two energy-conserving forms.
A relation analogous to Eq.~(\ref{eq:energy_semidiscr}) can be derived
by scalarly multiplying Eq.~(\ref{eq:NavStoDiscr2}) by $\mathbf{w}$, then by multiplying
the curl of Eq.~(\ref{eq:NavStoDiscr2}) by $\bu$, summing the two equations and using the symmetry of $\mathbf{R}$. The result is
\begin{equation}
\dfrac{1}{2}\dfrac{\text{d} h}{\text{d} t} = -\bu^T \mathbf{R} \bC(\bu) \bu 
- \cancel{\bu^T \mathbf{R}\mathbf{G} \bu} + \nu \bu^T \mathbf{R} \mathbf{L} \bu \;.
\end{equation}
The pressure term cancels as long as $\mathbf{R}\mathbf{G} = \mathbf{0}$, while the convective contribution depends again on the formulation employed for the convective term. In particular, the rotational form
preserves helicity, provided that the curl operator is symmetric. Indeed, upon substitution of the convective operator, $\bC(\bu) = -\bV( \bu)\mathbf{R}$, the nonlinear contribution reads
\begin{equation}
\bu^T \mathbf{R} \mathbf{V} (\bu) \mathbf{R} \bu = 0 \;,
\end{equation}
which cancels due to the skew-symmetry of $\mathbf{V}$ and the symmetry of $\mathbf{R}$.
On the other hand, the skew-symmetric, as well as the advective and divergence forms do not yield cancellation of the nonlinear term contribution in the helicity equation.
The above conclusions are in line with those that can be inferred by comparing Table~\ref{table:properties} and
Table~\ref{table:properties_discr}. Indeed, the skew-symmetric form requires the product rule to prove conservation of helicity, whereas for the rotational form only the symmetry of the curl operator is needed; the former property is not satisfied discretely (unless aliasing and truncation errors are removed), while the latter holds for any central-differencing method.

Although we do not attempt a formal proof, the fact that the discrete skew-symmetric form does not yield conservation of helicity, despite providing a skew-symmetric convective operator, can be intutively linked
to the nature of the helicity invariance. As mentioned in the Introduction, this is in fact not related to the skew-gradient nature of the Euler equations, but rather to the degeneracies of the Hamiltonian operator. In other words, both the skew-symmetric and the rotational forms yield a semi-discrete system that can be recast as Eq.~(\ref{eq:skewgrad}), but only one of them also possesses helicity as an inviscid invariant. 

\subsubsection{Staggered layout} \label{sec:staggered}

\begin{figure}
\centering
\includegraphics[width=.4\textwidth]{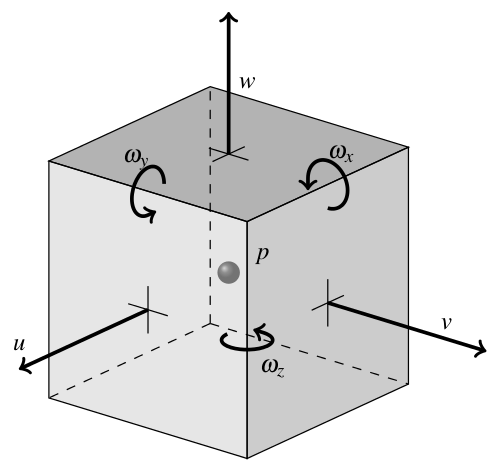}
\caption{Locations of the velocity components, vorticity components and pressure in a three-dimensional staggered layout. \label{fig:staggered}}
\end{figure}

The colocated arrangement is simple and straightforward to implement, although its
use is often complicated by the well-known odd-even decoupling phenomenon, which leads
to checkerboard patterns for the pressure field \cite{patankar1980numerical}. 
Nevertheless, the colocated layout is still
employed in many modern codes and special techniques have been developed to avoid 
the pressure-velocity decoupling \cite{rhie1983numerical,larsson2010co,trias2014symmetry}.
On the other hand, a common remedy is represented by the staggered layout, in which the pressure and the velocity
components are stored in different locations of the computational cell.
The most popular staggered arrangement was proposed by Harlow and Welch \cite{HW1965},
and is outlined in Fig.~\ref{fig:staggered}. 

A formal matrix analysis, similar to the one performed for the colocated case, 
becomes particularly involved for the staggered layout and is
beyond the scope of this paper. 
%For a demonstration of the energy-conservation properties
%of the Harlow-Welch scheme, the reader is referred to (REF).
Here, we focus on the scheme proposed by Harlow and Welch (HW) \cite{HW1965}, and test its
properties numerically in Section~\ref{sec:results}.
The HW scheme is a second-order method which relies on the divergence form of the convective term, as well as on proper interpolation operators; using the same notation as Morinishi \emph{et al.} \cite{Morinishi1998},
the nonlinear term is computed as
\begin{equation}
\text{(Div.)}_i \equiv \dfrac{\delta_1\overline{u}_j^{1x_i} \overline{u}_i^{1x_j}}{\delta_1 x_j} \;, \label{eq:HW}
\end{equation}
where $\delta_1$ and $\overline{\cdot}^{1x}$ are the half-grid derivative and interpolation
operators respectively. 

The HW scheme has several interesting properties: it preserves momentum a priori (i.e., regardless of mass conservation) as well as kinetic energy in the inviscid case. Indeed, it can be shown that the discretization in Eq.~(\ref{eq:HW}) leads to a skew-symmetric convective operator $\bC(\bu)$, provided that $\mathbf{M}\bu = \mathbf{0}$ \cite{coppolaAimeta}. On the other hand, to the authors' knowledge the helicity-conservation properties of the HW scheme have not been documented thus far in literature.
Note that the definition of helicity on a staggered grid is ambiguous, due to the relative positioning of velocity and vorticity components, see Fig.~\ref{fig:staggered}.
Here, we consider two possibilities, namely to compute helicity in the vertices and in the center of 
the computational cell
\begin{eqnarray}
    \left. h \right|_\text{vertices} = \overline{\overline{u}_1^{1x_2}}^{1x_3}\overline{\omega_1}^{1x_1} + \overline{\overline{u_1}^{1x_1}}^{1x_3}\overline{\omega_2}^{1x_2} + 
    \overline{\overline{u_3}^{1x_1}}^{1x_2}\overline{\omega_3}^{1x_3} \;, \label{eq:hel_stagg1}\\
    \left. h \right|_\text{centers} = 
    \overline{u_1}^{1x_1}\overline{\overline{\omega_1}^{1x_2}}^{1x_3} 
    + \overline{u_2}^{1x_2}\overline{\overline{\omega_2}^{1x_1}}^{1x_3}
    + \overline{u_3}^{1x_3}\overline{\overline{\omega_3}^{1x_1}}^{1x_2} \;. \label{eq:hel_stagg2}
\end{eqnarray}
We explicitly note that three interpolations have to be performed for each term to 
gather all the required variables in the same location. On the other hand,
among the various possible definitions, a form of kinetic energy that is preserved
by the HW scheme in absence of viscosity is the following \cite{ham2002fully}
\begin{equation}
    e = \dfrac{1}{2} \overline{u_i u_i}^{1x_i} \;, \\
\end{equation}
and is located at the cell center.

\subsection{Temporal discretization}

Once the Navier-Stokes equations have been spatially discretized, 
one is left with Eqs.~(\ref{eq:NavStoDiscr1})-(\ref{eq:NavStoDiscr2}), for  
which a further discretization in time is required.
The semi-discrete equations constitute an index-2 Differential Algebraic 
system \cite{Sanderse2012}, due to the presence of the pressure,
which acts as a kinematic constraint to ensure incompressibility.
A system of ODE is formally obtained upon introducing
a projection operator $\mathbf{P}$,
\begin{equation}\label{eq:ODESys}
\dfrac{\text{d} \bu}{\text{d} t}=\widetilde{\mathbf{f}}(\bu) \bu \;,
\end{equation}
where
$\widetilde{\mathbf{f}} = \mathbf{P}\mathbf{f}$ and
$\mathbf{f}\left(\mathbf{u}\right) = -\mathbf{C}\left(\mathbf{u}\right)+
\nu\mathbf{L}$,
with $\mathbf{P} = \mathbf{I}-\mathbf{G}\mathbf{\boldsymbol{\mathcal{L}}}^{-1}\mathbf{M}$ and
$\mathbf{\boldsymbol{\mathcal{L}}}=\mathbf{M}\mathbf{G}$.
Time advancement of Eq.~(\ref{eq:ODESys}) is now straightforward 
and can be carried out by means of any ODE solver. In general, the
discretization in time will in turn spuriously contribute to the
energy and helicity balance, unless proper care is taken.
Here we will focus on the class of Runge-Kutta (RK) time-advancement methods.
RK schemes are very popular in the fluid dynamics community due to their favorable 
properties, such as their self-starting capability and relatively large stability limit.

%The majority of turbulence simulations are nowadays performed by using
%three-stage (particularly the low-storage Wray's scheme \cite{Orlandi2000}) or four-stage (the
%classical RK4 \cite{GriffithsHigham}) methods in conjunction with fractional-step procedures.

A general $s-$stage Runge-Kutta method applied to Eq.~(\ref{eq:ODESys})
can be expressed as
\begin{eqnarray}
\bunp&=& \bun + \Delta t \sum_{i=1}^{s} b_{i} \widetilde{\mathbf{f}}(\bui)\bui \;, \label{eq:RK1} \\
%+ \Delta t \mathbf{G} \mathbf{p}^{n+1},            
\bui &=& \bun + \Delta t \sum_{j=1}^{s} a_{ij}
 \widetilde{\mathbf{f}}(\buj)\buj \;, \label{eq:RK2} %- \mathbf{G} \mathbf{p}_j  
\end{eqnarray}
where $a_{ij}$ and $b_i$ are the RK coefficients.
The RK coefficients are often arranged into the so-called
Butcher tableau \cite{Butcher2004}, and are usually constructed to maximize
the temporal order of accuracy of the method $p$; for $s \le 4$, schemes with
$p = s$ can be obtained.
Note that the role of the projection operator is to enforce incompressibility 
at each time instant and at each sub-step through the solution of a Poisson equation for pressure. This is often done in practice through a fractional-step method \cite{capuano2016approximate}.

Similarly to the spatial case, we first proceed to check whether the integration-by-parts
rule for the time-derivative operator, Eq.~(\ref{eq:SBPtime_cont}), is
retained by the RK method.
The discrete analogous of Eq.~(\ref{eq:SBPtime_cont}) reads
\begin{equation}
  \dfrac{\Delta (\hat{f} \cdot \hat{g})}{\Delta t} = \langle  \hat{f} \cdot \hat{G}  \rangle_{\Delta t} + \langle \hat{g} \cdot \hat{F} \rangle_{\Delta t} \;, \label{eq:SBPtime_discr}
\end{equation}
where $\hat{\cdot}$ indicates the time-discrete variable. In Runge-Kutta methods, the inner norm is properly defined as
\begin{equation}
	\langle f G \rangle_{\Delta t} \quad \rightarrow \quad \sum_{i=1}^s b_i  f_i G_i \;, \label{eq:innerNormTemp}
\end{equation}
which is reminiscent of the quadrature rule at the base of RK methods. 
Upon simple manipulations, one has
\begin{equation}
 \dfrac{\Delta (\hat{f} \cdot \hat{g})}{\Delta t} =  \langle  \hat{f} \cdot \hat{G}  \rangle_{\Delta t} + \langle \hat{g} \cdot \hat{F} \rangle_{\Delta t} - 
 \Delta t \sum_{i,j} \left(b_ia_{ij}+b_j a_{ji}-b_ib_j\right) \hat{F}_i \hat{G}_j \;. \label{eq:RK_symp}
\end{equation}
In practice, upon comparison of Eq.~(\ref{eq:RK_symp}) with Eq.~(\ref{eq:SBPtime_discr}), it is evident that RK schemes do not satisfy a discrete integration-by-parts rule unless the identity
\begin{equation}
g_{ij} \equiv b_ia_{ij}+b_j a_{ji}-b_ib_j = 0 \label{eq:cond_symp}
\end{equation}
is verified for every $i,j$. The condition (\ref{eq:cond_symp}) is very
well known in the community of so-called geometric integration methods,
and is the requirement for a RK method to preserve all the quadratic invariants
of a system \cite{Hairer2006}. Also, for irreducible RK methods, Eq.~(\ref{eq:cond_symp})
is a necessary and sufficient condition for the method to be \textit{symplectic}.
Schemes satisfying Eq.~(\ref{eq:cond_symp}) are necessarily implicit.

%where $f_i$ and $G_i$ follow from Eq.~(\ref{eq:RK2}). In non-symplectic methods, Eq.~(\ref{eq:SBPtime_discr}) does not hold.

The derivation of fully-discrete balance equations for energy and helicity
is now straightforward. Again, the expression for the kinetic energy variation
introduced by Eqs.~(\ref{eq:RK1})-(\ref{eq:RK2}) has been derived 
elsewhere \cite{Capuano2015b,Sanzserna,Sanderse2013} and is reported here for convenience,
\begin{equation}
\dfrac{\Delta e}{\Delta t} = 
\nu \sum_i^s b_i \bu^T_i \mathbf{L} \bu_i 
- \dfrac{\Delta t}{2} \sum_{i,j}^s g_{ij} \bu_i^T
 \widetilde{\mathbf{f}}^T\left(\bui\right)\widetilde{\mathbf{f}}\left(\bu_j\right) \bu_j \;,
 \label{eq:energy_fullydiscr}
\end{equation}
where it has been assumed that the spatial discretization is energy conserving.

We focus hereinafter on the fully discrete evolution equation for helicity,
which reads
\begin{equation}
%\begin{split}
\frac{\Delta h}{\Delta t} = 2\nu \sum_i^s b_i \bu_i^T \mathbf{R}  \mathbf{L} \bu_i
 -2 \sum_i^s b_i \bu_i^T \mathbf{R} \mathbf{C} \left(\bu_i \right) \bu_i 
 -\Delta t\sum_{i,j}^s g_{ij}  \widetilde{\mathbf{f}}_i^T  \mathbf{R}  \widetilde{\mathbf{f}}_j \;.
\label{eq:helicity_fullydiscr}
%\end{split}
\end{equation}
Note that the equations of both energy and helicity resemble the generic structure given in Eq.~(\ref{eq:RK_symp}).
Particularly, the three terms appearing in the right-hand side of Eq.~(\ref{eq:helicity_fullydiscr}) are respectively the contribution of the
discretized (physical) viscous dissipation, the spatial error term due to convection, and the temporal error. The spatial error behaviour follows from the results already obtained in Section~\ref{sec:spatial}, e.g., it vanishes for the rotational form on a colocated layout. The temporal error
stems from the possible lack of the integration-by-parts rule. Therefore, only symplectic RK methods are able to preserve energy and helicity in time in case of inviscid flow. This is not surprising, since these two quantities are inviscid invariants of the Navier-Stokes equations, and, as mentioned above, RK methods satisfying $g_{ij} = 0$ preserve all the quadratic invariants of a system.
When energy- and helicity-preserving methods are employed both in space and time, Eqs.~(\ref{eq:energy_fullydiscr})-(\ref{eq:helicity_fullydiscr}) become the discrete counterparts of Eqs.~(\ref{eq:energy_int})-(\ref{eq:helicity_int}), i.e., the induced balance equations of energy and helicity are correctly enforced, regardless of scheme order, grid size and time step \cite{Sanderse2013}. On the other hand, the discrete dissipation rate is obviously still affected by spatial and temporal discretization errors, when compared to its exact counterpart.

As far as time-integration is concerned, it is worth to note that recently Capuano \emph{et al.} developed special explicit Runge-Kutta methods for the fluid flow equations, named \textit{pseudo-symplectic}, with the aim of preserving kinetic energy in time up to an order $q$, with $q > p$. Nonetheless, these schemes have actually order $q$ for the conservation of \textit{any} quadratic invariant, and are therefore applicable also for the purpose of preserving helicity (or enforcing the correct helicity balance) in time to a higher order of accuracy.
%\textcolor{red}{come abbiamo mostrato nelle nostre simulazioni, gli errori di discretizzazione modificano comunque le u, anche se la discretizzazione conserva energia/elicita' globali. percio' il campo modificato puo' dissipare piu o meno rispetto ad una soluzione esatta, e quindi ci sono variazioni di energia/elicita globali! anzi questo e'  uno dei nostri messaggi. io eliminerei da 'When energy...' in poi e scriverei qlks del tipo: When energy- and helicity-preserving methods are employed both in space and time in an inviscid simulation, the induced balance equations of energy and helicity are satisfied to machine precision, regardless of scheme order, grid size and time step. However, it is worth noticing that the discretization error still affects the velocity unknowns, and therefore in the viscous case (where the velocity-dependent viscous term in Eqs.~\eqref{eq:energy_fullydiscr}-\eqref{eq:helicity_fullydiscr} does not vanish) even a fully conservative method displays a numerical error in terms of global energy and helicity. This will be further investigated through a scale-dependent analysis of discretization error in section ...}

By collecting the outcome of the preceding sections with known literature results \cite{Morinishi1998,Capuano2015b}, the conservation properties (also including mean momentum and energy) of
the Navier-Stokes algorithms analyzed in the paper are reported in Table~\ref{table:consprop}.
Any explicit RK method produces an error in the conservation of quadratic invariants, although the mean momentum is preserved, since all Runge-Kutta methods preserve linear invariants \cite{Rosenbaum1976}. Only the rotational and the skew-symmetric forms on a colocated layout, and the staggered HW scheme are able to preserve kinetic energy,
although the condition $\mathbf{M}\bu = \mathbf{0}$ (as well as the usual consistency relation $\mathbf{G}^T = -\mathbf{M}$) must hold.
Remarkably, both the convective and the pressure terms in the helicity equation vanish regardless of the discrete enforcement of continuity, and therefore a discretization based on the rotational form preserves helicity a priori (assuming $\mathbf{R}\mathbf{G} = \mathbf{0}$). Note that the rotational form, coupled to a symplectic integrator, is able to preserve mean momentum, energy and helicity in the inviscid case.

\begin{table}
\centering
  \caption{Conservation properties of mean momentum, energy and helicity for 3D Navier-Stokes discretizations, for various forms of the nonlinear term and symplectic (symp.) or explicit (expl.) Runge-Kutta schemes. Algorithms 1-5 refer to a colocated grid, Algorithm 6 is the Harlow-Welch scheme described in Section~\ref{sec:staggered}. 
  The discrete integration-by-parts rule is assumed to hold.
  %for global momentum $M$, energy $E$ and helicity $H$. S: space, T: time. 
  $+$ conservative \textit{a priori}, $\circ$ conservative \textit{if and only if} continuity is discretely satisfied, $\times$ non conservative.
   \label{table:consprop}}
  \begin{tabular}{lllllllll}%{p{1.3cm}p{1.5cm}p{1.5cm}p{1cm}p{1cm}p{1cm}p{1cm}p{1cm}p{1cm}}
   \hline\noalign{\smallskip}
    & \multicolumn{2}{l}{Algorithms} & \multicolumn{2}{p{2cm}}{Momentum} & \multicolumn{2}{p{2cm}}{Energy} & \multicolumn{2}{p{2cm}}{Helicity} \\
   %\cmidrule{2-9}
   \# & Space & Time & Space & Time & Space & Time & Space & Time \\
   \noalign{\smallskip}\hline\noalign{\smallskip}
   1 & Rot. & Symp. & $\circ$ & $+$ & $\circ$  & $+$    & $+$      & $+$ \\
   2 & Rot. & Expl. & $\circ$ & $+$ & $\circ$  & $\times$ & $+$      & $\times$  \\
   3 & Skew & Symp. & $\circ$ & $+$ & $\circ$  & $+$      & $\times$ & $+$ \\
   4 & Skew & Expl. & $\circ$ & $+$ & $\circ$  & $\times$ & $\times$ & $\times$  \\
   5 & Div. & Expl. & $+$     & $+$ & $\times$ & $\times$ & $\times$ & $\times$ \\
   6 & HW   & Expl. & $+$     & $+$ & $\circ$ &  $\times$ & $\times$ & $\times$ \\
   \noalign{\smallskip}\hline
  \end{tabular}
\end{table}

\subsection{Extension to nonuniform meshes} \label{sec:irregular}

The hypothesis of uniform Cartesian grids made in the preceding sections
can be easily removed by taking into account a relevant inner product,
\begin{equation}
\langle uv \rangle \quad \rightarrow \quad \bu^T \bO \mathbf{v} \;,
\end{equation}
where $\bO$ is a diagonal matrix containing the metrics of the mesh.
In this case Eq.~(\ref{eq:NavStoDiscr2}) is most properly recast as
follows
\begin{equation}
\bO \frac{\text{d} \bu}{\text{d} t} + \bC(\bu)\bu = - \bO \mathbf{G}\mathbf{p}
+ \nu \mathbf{L}\bu \;, \label{eq:NavStoDiscr2Irr}
\end{equation}
and the requirements for energy and helicity conservation can
be easily shown to be the same as in the uniform case, i.e., a skew-symmetric convective operator $\bC(\bu)$ for the energy, and a symmetric curl operator $\mathbf{R}$
for the helicity. The proper relation between the gradient and the divergence operators in this case is $\mathbf{G}^T \bO = -\mathbf{M}$, while $\mathbf{R}\mathbf{G} = \mathbf{0}$ must still hold. With a careful choice of the discrete differential operators, 
extension to general unstructured meshes is also possible in principle, but will not be
detailed here.

%In case of nonuniform, non-Cartesian grid, care must be 
%taken in constructing the operators so that the aforementioned properties
%are enforced. However, 

\section{Numerical results} \label{sec:results}

All the numerical results presented in this section have been produced by a pseudo-spectral code in a three-dimensional periodic box. This choice has allowed very efficient computations and the possibility to analyze the numerical errors in a scale-dependent way through spectra of energy, helicity and their non-linear transfers. It was also possible to analyze the effects of truncation errors by mimicking standard non-dealiased finite-difference schemes, on both a colocated and a staggered arrangement, via the modified wavenumber approach, and to produce accompanying de-aliased spectrally resolved computations for comparison. It is worth to recall that the modified wavenumber approach consists in substituting the wavenumber vectors used for the computation of spectrally-accurate derivatives with their finite-difference counterparts \cite{Kravchenko1997}. For instance, for a second-order colocated scheme, the modified wavenumbers for the first- and second-derivative read respectively
\begin{equation}
k'(k) = \dfrac{\sin \left( k h \right) }{h} \;, \quad \quad 
k'^2(k) = \dfrac{2 \left(1-\cos\left(k h \right) \right)}{h^2} \;,
\end{equation}
where $h$ is the grid size.
The staggered scheme has been mimicked by introducing proper shift and averaging operators to
take into account the necessary interpolation procedures.
For all the simulations the time step $\Delta t$ is selected based on a CFL-like condition with CFL=0.5.

\subsection{Inviscid helical dynamics}

As a first test, we have simulated an inviscid system with the aim of isolating the errors coming from the discretization of the nonlinear term and the time-advacing scheme.
We have thus solved the Euler equations on a grid composed by $32^3$
points, and tested the various algorithms listed in Table~\ref{table:consprop}
coupled with various Runge-Kutta methods and both with second-order and spectral accuracy.
The initial condition is the superposition of two Arnold-Beltrami-Childress (ABC) flows \cite{childress2008stretch} at the wavenumbers $k_1=4$ and $k_2 = 6$.

\begin{figure}
\centerline{
 \includegraphics[height=0.42\textwidth]{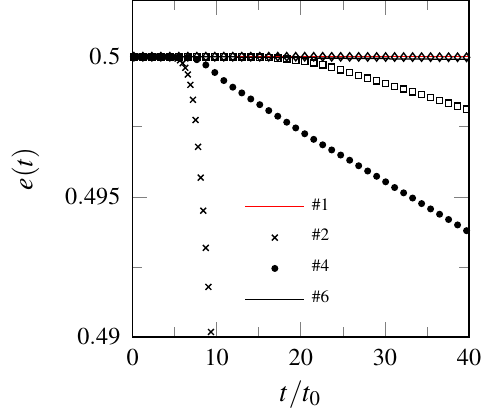}
 \hfill
 \includegraphics[height=0.42\textwidth]{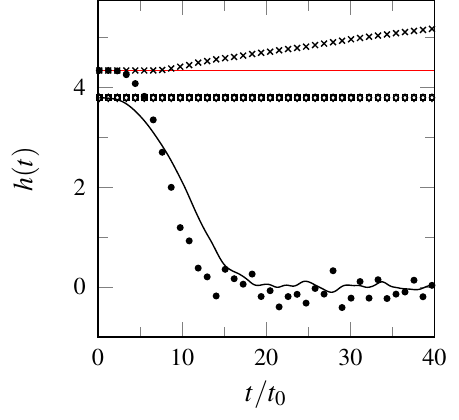} % using a pdf figure and PdfTex
 } % using a pdf figure and PdfTex
\caption{Time evolution of global energy and helicity for Euler dynamics using the various algorithms listed in Table~\ref{table:consprop}. In all cases, spectral accuracy and a RK4 scheme are used, except for two variants of the Algorithm 2, namely with II order accuracy (square) and with II order accuracy and the 4p7q RK scheme \cite{capuano2017explicit} (diamond). \mbox{$t_0= e^{-1/2} k_1^{-1}$} is the system characteristic time scale.}
\label{fig:hw_evolutions}
\end{figure}

The time evolutions of energy and helicity are shown in the left and right part of Fig.~\ref{fig:hw_evolutions} respectively. The colocated algorithms have been used with full spectral accuracy and with the classical fourth-order Runge-Kutta scheme (RK4), except for two runs, in which the rotational form has been coupled to a second-order spatial scheme, and to both a second-order scheme and a pseudo-symplectic RK method. Particularly, the explicit six-stage scheme with $p=4$ and $q=7$ (named 4p7q) has been tested; the reader is referred to \cite{capuano2017explicit} for details. The symplectic scheme in Algorithm \#1 is the Gauss midpoint method. 

The results fully confirm the theoretical predictions.
The rotational form preserves helicity spatially, although it produces significant deviations in both quadratic invariants when used with an explicit time-advancement scheme. This behaviour is attributed to the prominent accumulation of energy at smaller scales typical of spectrally truncated conservative inviscid systems, which causes the temporal error of the RK method to be significant. Remarkably, the RK error is dissipative for energy (as it is well known), but productive for helicity. A possible explanation for this circumstance might be the fact that the temporal error in the discrete helicity equation does not have a definite sign and therefore helicity can increase. When the rotational form is used in conjunction with a second-order scheme, the errors are much less pronounced, likely due to the attenuation of the high-wavenumber energy content by the finite-difference derivative method. In this regard, note that the group of second-order runs has a different initial helicity content due to truncation errors affecting the computation of initial vorticity. The beneficial impact of the higher-order conservation properties of the pseudo-symplectic RK method is evident in keeping the two invariants practically constant, as done to machine precision by the fully conservative Algorithm 1.

The skew-symmetric form is shown to dissipate the initial helicity content in few characteristic times. On the other hand, on equal spatial accuracy, the skew-symmetric form has a better behaviour on energy than the rotational form, due to its well-known favorable aliasing cancellation properties \cite{Kravchenko1997}. When coupled to a lower-order scheme (not shown here), the behaviour of the skew-symmetric form is qualitatively very similar to the one already described for the rotational formulation. 

The staggered Harlow-Welch scheme is found to dissipate helicity, with a rate similar to the colocated skew-symmetric form. Actually, both the expressions proposed to compute helicity in a staggered layout, Eqs.~\eqref{eq:hel_stagg1}-\eqref{eq:hel_stagg2} have been found to be practically coincident from a numerical point of view. It can thus be concluded that the classical energy-conserving second-order staggered method does not preserve helicity, at least with the definitions of vorticity and helicity given in Section~\ref{sec:spatial}. The colocated divergence form (Algorithm 5) diverged immediately due to violation of energy conservation.

The analysis reported above has allowed to confirm the theoretical results and gain preliminary insights into the behaviour of the two energy-conserving forms. With the aim of obtaining scale-dependent information, we have computed energy and helicity spectra for the solution of spherically truncated Euler equations, i.e. the system obtained by truncating the Euler equations in Fourier space at a wavenumber \sloppy{$|\boldsymbol{k}| = k_\textrm{\scriptsize max}$}. In the absence of numerical errors, this truncated system preserves energy and helicity. In correspondence of the statistically stationary state (\textit{absolute} or statistical equilibrium), exact expressions for the spectra of both energy and helicity are available \cite{kraichnan1973helical},
\begin{equation}
  E(k)=\frac{4\pi}{\alpha}\frac{k^2}{1-\big(\frac{\beta}{\alpha}\big)^2k^2}\;,
  \quad \quad
  H(k)=\frac{8\pi\beta}{\alpha^2}\frac{k^4}{1-\big(\frac{\beta}{\alpha}\big)^2k^2}\;,
\label{eq:exactspectra}
\end{equation}
where $\alpha$ and $\beta$ are constants that depend on the energy and helicity contents.
The statistical equilibrium is governed by a single non-dimensional parameter, e.g. the relative helicity $H_\textrm{\scriptsize rel}=K_h/(2 k_\textrm{\scriptsize max}K_e)$.
In the non-helical case, i.e. $H_\textrm{\scriptsize rel}=\beta=0$, the helicity spectrum vanishes and the energy spectrum becomes simply
\begin{equation}
  E(k)=\frac{4\pi}{\alpha}k^2 \;.
\label{eq:E(k)}
\end{equation}
In the performed simulations $k_\textrm{\scriptsize max}=42$ and the initial condition is again the sum of two Arnold-Beltrami-Childress (ABC) flows with $k_1 = 28$ and $k_2 = 30$, similarly to the test proposed in \cite{krstulovic2009cascades}, 
yielding $H_\textrm{\scriptsize rel}=0.687$.
The same algorithms described in the previous case are tested here with exception of the case with the pseudo-symplectic method.

%Time evolutions of global energy and helicity are reported in Fig.~\ref{fig:evolutions}. Algorithm 1 (rotational form, symplectic RK schemes) is fully conservative and indeed is found to preserve energy and helicity in time up to machine precision. The rotational and skew-symmetric forms in conjunction with explicit Runge-Kutta schemes (Algorithms 2 and 4) are found to slightly dissipate energy in time, with Algorithm 2 being more dissipative than Algorithm 4. This is attributed to the accumulation of energy at the smallest scales due to the larger aliasing errors of the rotational form.
%On the other hand, the skew-symmetric form is found to completely dissipate the initial helicity contents in few characteristic times. It is interesting to note that Algorithm 2, which employs the rotational form and an explicit RK3, is slightly helicity-productive. We conjecture that this might be due to the fact that the temporal error in the discrete helicity equation does not have a definite sign and thus might be productive.

%\begin{figure}
%\centerline{
% \includegraphics[height=0.42\textwidth]{figure1.eps}
% \hfill
% \includegraphics[height=0.42\textwidth]{figure2.eps} % using a pdf figure and PdfTex
% } % using a pdf figure and PdfTex
%\caption{Time evolution of global energy and helicity for spherically truncated Euler dynamics using the algorithms listed in Table~\ref{table:consprop}; $t_0= K_e^{-1/2} k_\textrm{\scriptsize max}^{-1}$ is the system characteristic time scale.}
%\label{fig:evolutions}
%\end{figure}

\begin{figure}
\centerline{
 \includegraphics[height=0.42\textwidth]{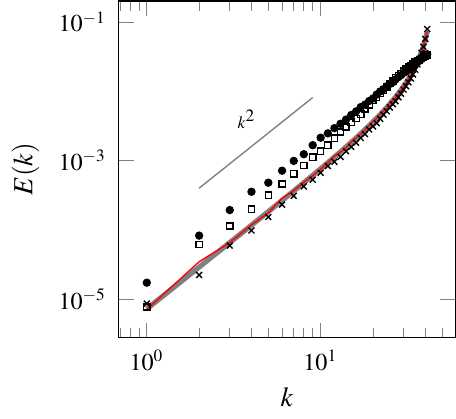}
 \hfill
 \includegraphics[height=0.42\textwidth]{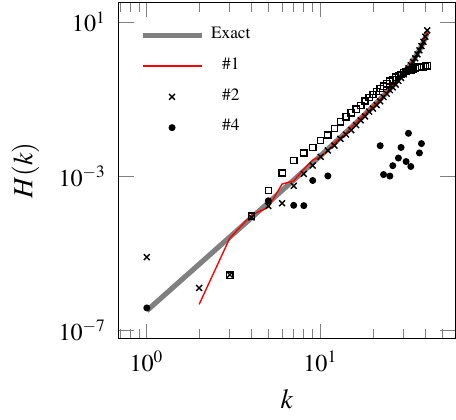} % using a pdf figure and PdfTex
 } % using a pdf figure and PdfTex
\caption{Energy and helicity spectra in spherically truncated Euler dynamics for various algorithms, see Table~\ref{table:consprop} and the label of Figure~\ref{fig:staggered} for further details. The square symbol refers to Algorithm 2 with second-order spatial accuracy. The grey curves display the exact solution which is 
reported in Eq.~(\ref{eq:exactspectra}). }
\label{fig:spectra}
\end{figure}

Energy and helicity spectra are reported in Fig.~\ref{fig:spectra}. The fully conservative computation matches very closely the exact solution, while Algorithm 2 (rotational form with spectral accuracy) is slightly dissipative due to the temporal dissipation of the explicit scheme. Importantly, the computation employing Algorithm 4 (skew-symmetric form with spectral accuracy, explicit) drives the system towards the non-helical equilibrium solution, i.e. the relative helicity vanishes and $E(k)$ is proportional to $k^2$, see Eq.~(\ref{eq:E(k)}). 
The behaviour of Algorithm 2 when the accuracy of the spatial derivatives is reduced to second order deserves attention. Although this method preserves energy and helicity spatially and with overall good accuracy (even better than its spectral version, see Fig.~\ref{fig:staggered}), the energy and helicity spectra are distorted with respect to the analytical solution, and present significant errors over almost all the wavenumber range, which can be attributed to the interaction of the truncation error with the mechanism of nonlinear helicity transfer.
We conjecture that this distorsion is also responsible for the aforementioned reduced energy and helicity dissipation arising from the temporal error.
However, a cross-comparison with a low-order scheme based on the skew-symmetric form cannot be performed in this inviscid framework since the helicity is quickly dissipated.

\subsection{Forced helical turbulence}

In this section, the effect of the spatial discretization errors related to the convective term is investigated through viscous helical simulations. In all cases, time-advancement is performed through a fourth-order RK scheme.

We choose to force our simulations at large scales in order to achieve a statistically stationary state and then compute spectra by time averaging. 
The selected forcing scheme is the helical Euler forcing \cite{kerr_1985,vallefuoco2017}, which is inspired by the truncated Euler dynamics: the lowest-wavenumber modes, corresponding to wavevectors $\boldsymbol{k}$ such that $|\boldsymbol{k}| \le k_F$ ($k_F$ is the largest forcing wavenumber), obey the three-dimensional incompressible Euler equations and are independent of the other modes.
Therefore, we set $k_F=2.5$ and inject both energy and positive helicity in the range $k\le k_F=2.5$. The relative helicity of the Euler system, defined as $\Hc^E/(2 k_F \Ec^E)$, is equal to \mbox{$H_\text{rel}=0.967$}.

Energy and helicity are then transfered to higher wavenumbers down to the dissipative scales, where they are dissipated.
It is worth to recall the Lin equations for energy and helicity spectra, which read
\begin{eqnarray}
\frac{\partial E(k)}{\partial t}=T_e(k)-2\nu k^2E(k)+\Phi_e(k) \;, \label{eq:LinE} \\
\frac{\partial H(k)}{\partial t}=T_h(k)-2\nu k^2H(k)+\Phi_h(k) \;,   \label{eq:LinH}
\end{eqnarray}
where $\Phi_e(k)$ and $\Phi_h(k)$ are linked to the external forcing and vanish in the unforced range $k > k_F=2.5$.
The non-linear transfers $T_e(k)$ and $T_h(k)$ represent the ratio of net energy and helicity transfered to scale $k$ through non-linear interactions among modes. 
For an energy- or helicity-conserving discretization, the sum of $T_e(k)$ or $T_h(k)$ over all scales vanishes, while a non conservative discretization leads to a net global transfer.
Under the assumption of statistically stationary flow, and restricting the balance to the unforced scales, Eqs.~\eqref{eq:LinE}-\eqref{eq:LinH} reduce to
\begin{eqnarray}
T_e(k)=2\nu k^2E(k) \;, \label{eq:LinEred} \\
T_h(k)=2\nu k^2H(k) \;,   \label{eq:LinHred}
\end{eqnarray}
that is, energy and helicity transfers are balanced by viscous dissipation.

Three simulations are performed in order to study the behaviour of the rotational and the skew-symmetric forms: a reference simulation, spectrally resolved and fully de-aliased (2/3 rule) with a 192$^\text{3}$ resolution (that is, 128$^\text{3}$ effective resolution); a second-order 128$^\text{3}$ simulation employing the rotational form; and a second-order 128$^\text{3}$ simulation employing the skew-symmetric form.
Table~\ref{table:helruns} reports some relevant parameters of these three runs.

\begin{table}
\centering
  \caption{Parameters of helical turbulence simulations. $\Ec^E$ and $\Hc^E$ are the energy and helicity contents of the Euler system. $\epsilon_e$ and $\epsilon_h$ are the energy and helicity average dissipation rates. 
  \label{table:helruns}}
  \begin{tabular}{lllll}
   \hline\noalign{\smallskip}
   {Algorithms} & {$e/\Ec^E$} & {$h/\Hc^E$} & $\sum T_e(k)/\epsilon_e$ & $\sum T_h(k)/\epsilon_h$  \\
   \noalign{\smallskip}\hline\noalign{\smallskip}
  Deal. & 2.16 & 1.82 & -0.012 & -0.021  \\
  Rot. & 2.22 & 1.95 & -0.015 & -0.0032   \\
  Skew. & 2.19 & 1.86 & -0.005 & -0.13 \\
   \noalign{\smallskip}\hline
  \end{tabular}
\end{table}

Spectra and statistical quantities are obtained by time-averaging over 310  turnover times after the statistically stationary state has been reached. 
The large-scale resolution is $L/(2 \pi)=0.11$ ($L$ is the integral scale and $2\pi$ is the periodic domain size), while the small scale resolution is $k_\text{max}\eta=1.12$ ($\eta$ is the Kolmogorov scale and $k_\text{max}=64$ is the maximum resolved wavenumber).
The Reynolds number is $\text{Re}=L\sqrt{2/3 e}/\nu=173$.

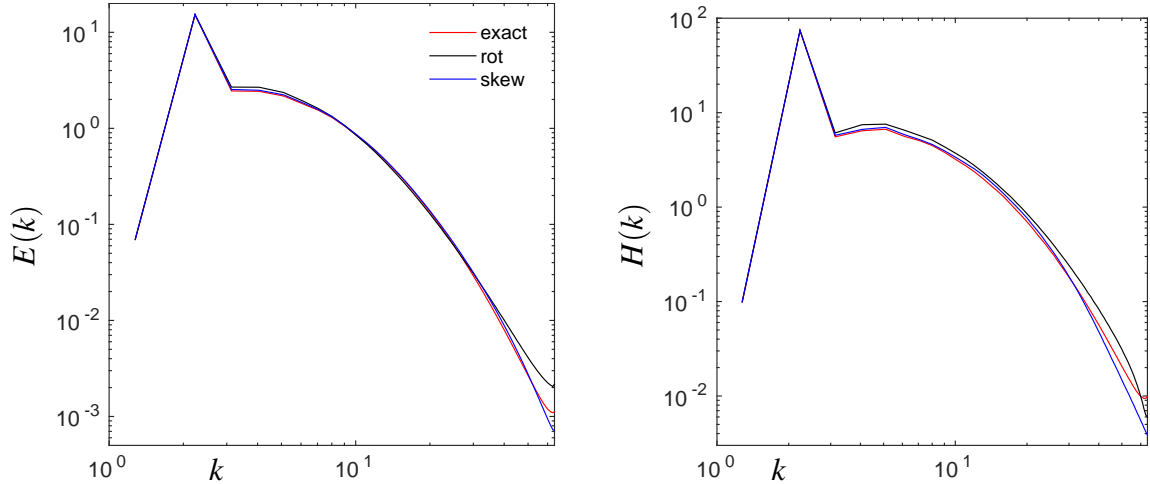
\begin{figure}
\begin{tikzpicture}[]
         \node at (-3.5,1.2) {\includegraphics[width=.53\textwidth]{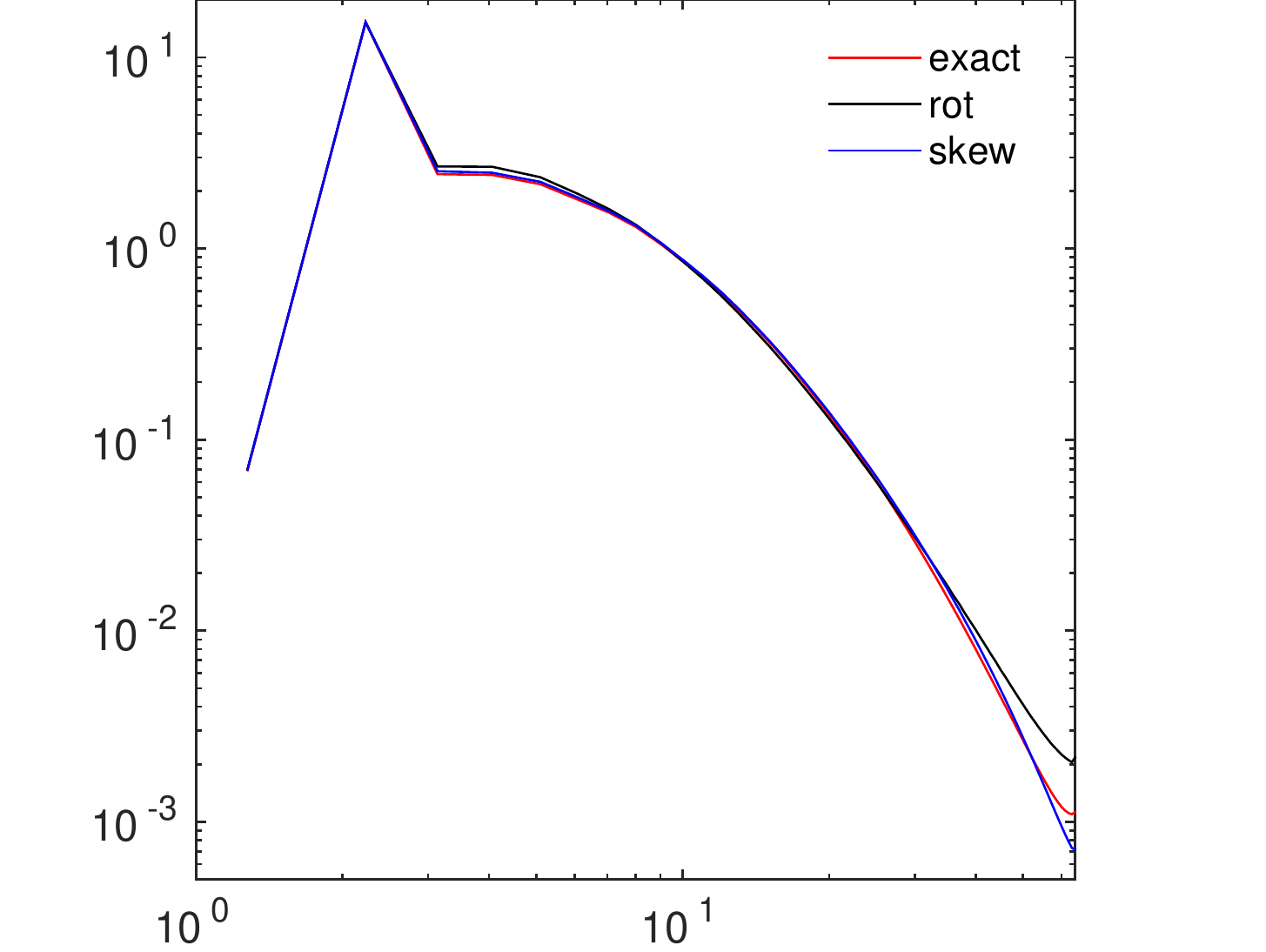}};
         \node at (4.4,1.2) {\includegraphics[width=.53\textwidth]{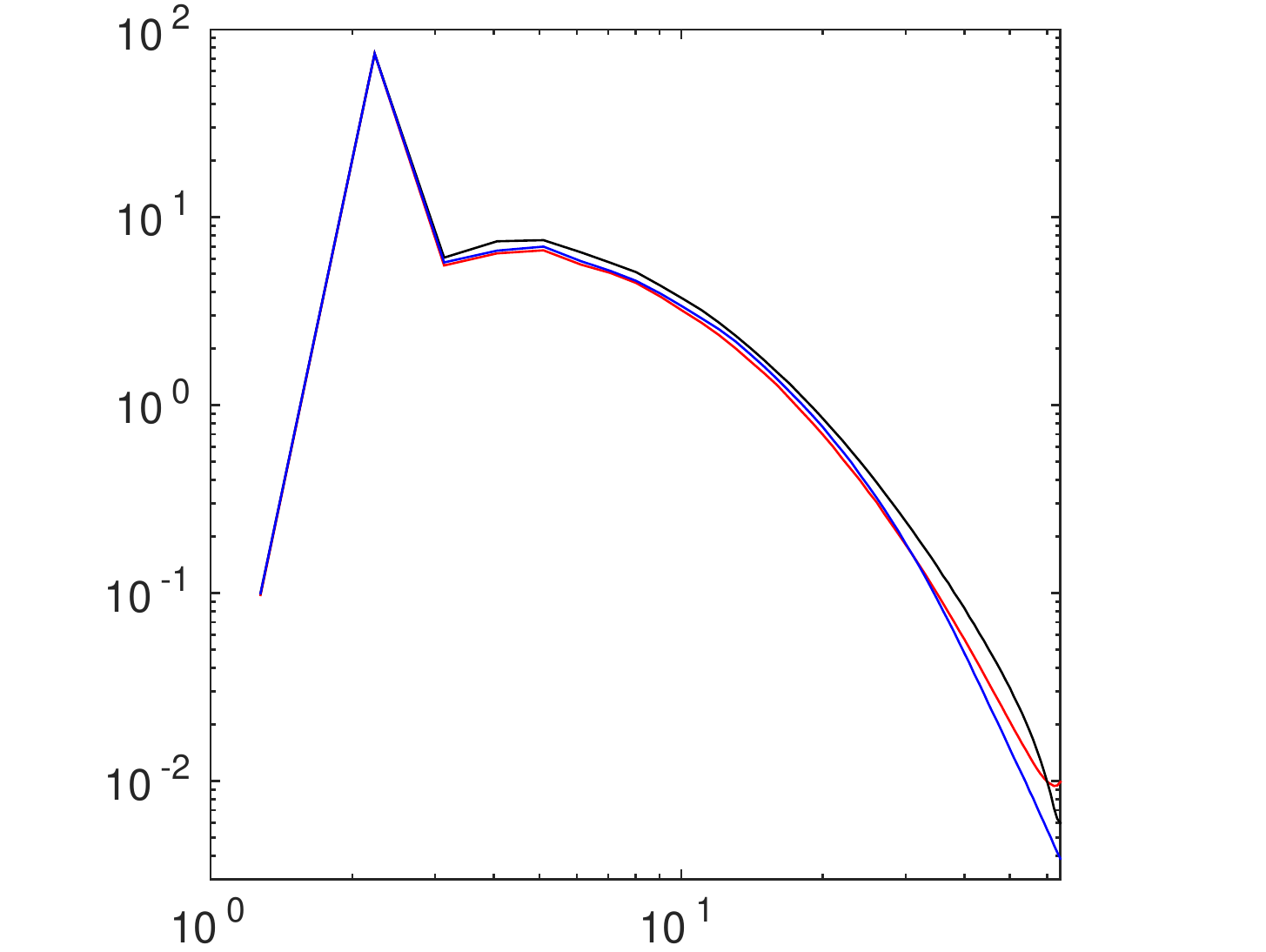}};
         \node at (-5,-1.8) {\large$k$};
         \node at (2.4,-1.8) {\large$k$};
         \node [rotate=90] at (-7.5,1.3) {\large$E(k)$};
         \node [rotate=90] at (.5,1.3) {\large$H(k)$};
\end{tikzpicture}
\caption{Energy and helicity spectra in helical turbulence for the de-aliased spectrally resolved formulation (red), the rotational formulation (black) and the skew-symmetric formulation (blue).\label{fig:EHspectra}}
\end{figure}

The energy and helicity spectra are plotted in Fig.~\ref{fig:EHspectra}.
The energy spectra of both the rotational and skew-symmetric algorithms almost superpose with the reference solution spectrum.
Only at the smallest resolved scales, in the dissipative range, the discretization error in the rotational simulation causes an energy pile-up, in accordance with previous observations \cite{Zang1991}. On the other hand, the departure of the skew-symmetric algorithm from the exact energy spectrum is significantly smaller, and of the same order as the slight pile-up of the reference solution.
The small departure from the exact solution in terms of energy is confirmed by the values of the global contents of energy reported in Table~\ref{table:helruns}, the relative error for the rotational form being twice the one yielded by the skew-symmetric.

As for the helicity spectra in Fig.~\ref{fig:EHspectra}, the differences among the three simulations are larger.
Surprisingly, the rotational form seems to yield a larger departure from the exact helicity spectrum with respect to the skew-symmetric form, despite its helicity-preserving property.
In particular, for the skew-symmetric form the helicity spectrum is negligibly larger than the reference one at low wavenumbers, while being significantly lower at smaller scales. On the other hand, for the rotational form, helicity is clearly too large over almost all the resolved range, becoming lower than the reference values only at the smallest resolved scales, due to the slight helicity pile-up in the de-aliased spectrally resolved simulation. 
These observations are confirmed by the global helicity contents reported in Table~\ref{table:helruns}.

At this point, it is worth to investigate why the discretization error in the (helicity-conserving) algorithm based on the rotational form affects the helicity spectrum more than for the skew-symmetric case, which does not conserve helicity.
In this regard, it is worth to preliminarily remind that although energy- and helicity-preserving methods satisfy the global balance equations to machine precision,
the discretization error still affects the velocity unknowns, and thus might alter both the scale-by-scale nonlinear transfer as well as the viscous dissipation rate.
In light of these clarifications, it appears from Fig.~\ref{fig:EHspectra} that, while discrete energy conservation leads to accurate energy spectra, the helicity-preserving property is not sufficient to reproduce the correct helicty spectrum. Therefore, we proceeded to carry out a scale-dependent conservation analysis to further characterize the effect of discretization errors on helicity.

%the important result arising from Fig.~\ref{fig:EHspectra} is that, while energy conservation leads to accurate energy spectra, the helicity-conserving property is not enough for a numerical scheme to reproduce the correct helicty spectrum. Therefore, a scale-dependent conservation analysis is required to characterize the effect of discretization errors on helicity.

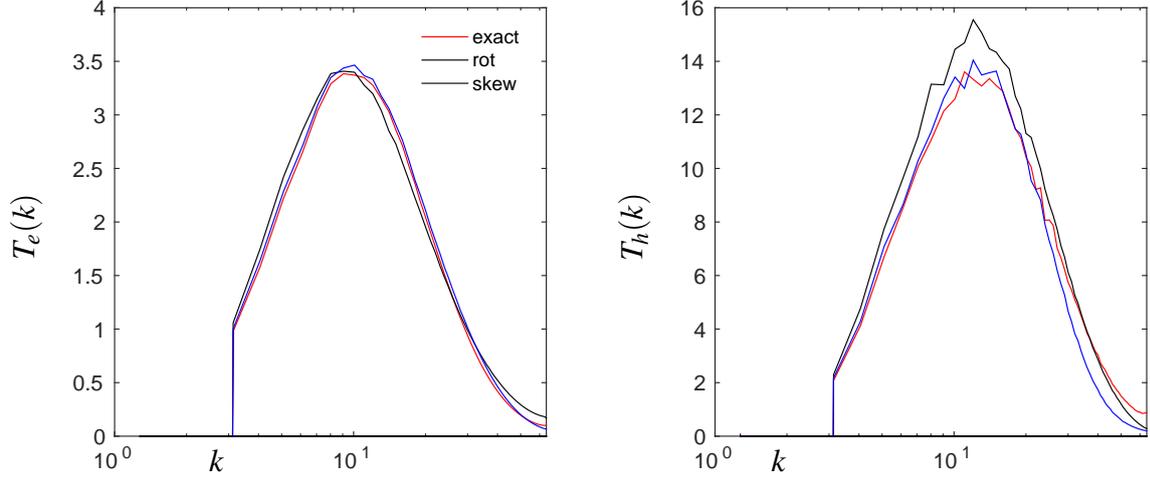
\begin{figure}
\begin{tikzpicture}[]
         \node at (-3.5,1.2) {\includegraphics[width=.53\textwidth]{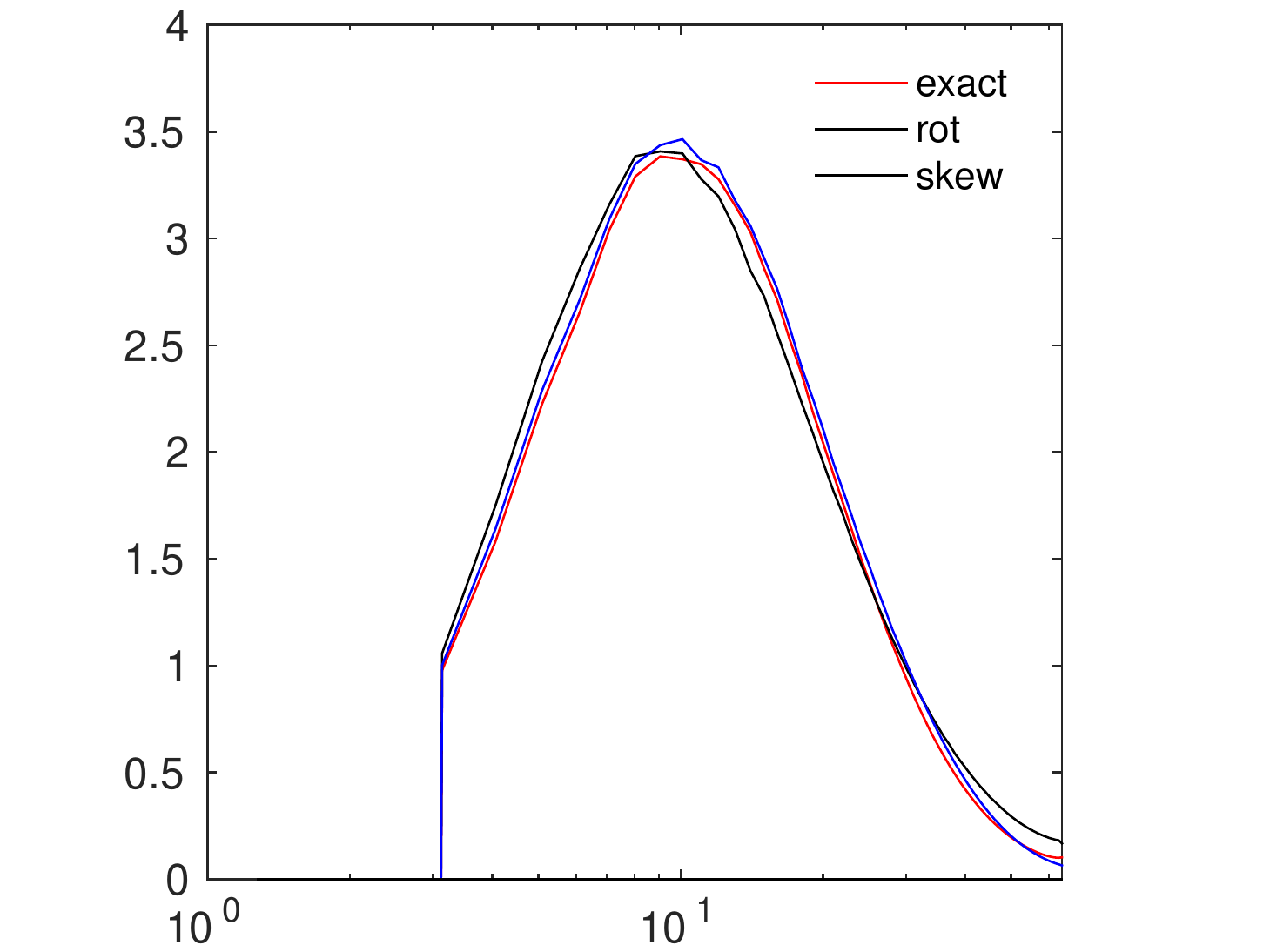}};
         \node at (4.4,1.2) {\includegraphics[width=.53\textwidth]{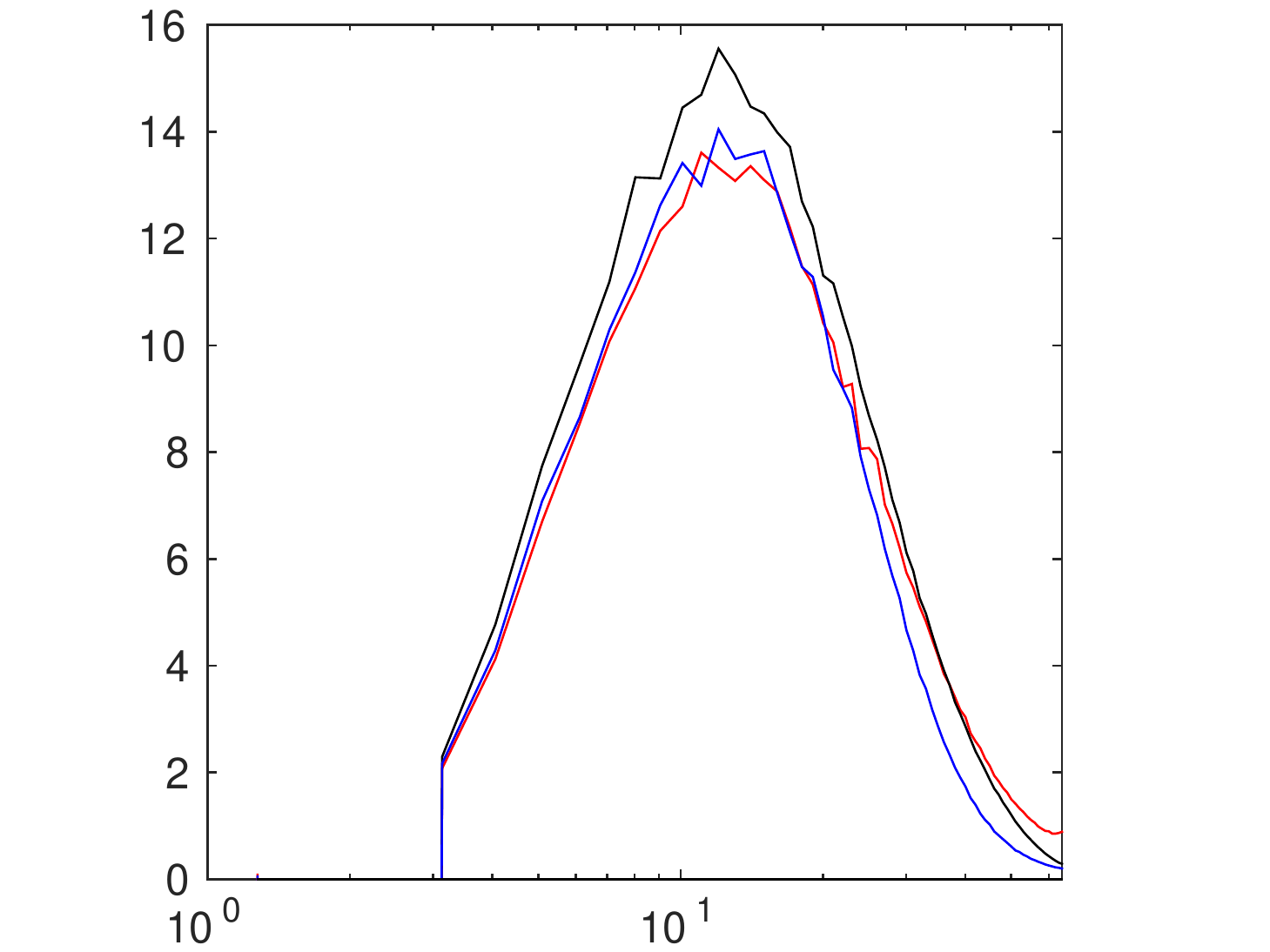}};
         \node at (-5,-1.8) {\large$k$};
         \node at (2.4,-1.8) {\large$k$};
         \node [rotate=90] at (-7.5,1.3) {\large$T_e(k)$};
         \node [rotate=90] at (.5,1.3) {\large$T_h(k)$};      
\end{tikzpicture}
\caption{Energy and helicity transfer spectra in helical turbulence for the de-aliased spectrally resolved formulation (red), the rotational formulation (black) and the skew-symmetric formulation (blue).\label{fig:Tspectra}}
\end{figure}

We thus compute the energy and helicity non-linear transfers $T_e(k)$ and $T_h(k)$.
Table~\ref{table:helruns} reports the sums of the transfers over scales $0\le k \le k_\text{max}$. For the energy and helicity transfers in the spectrally resolved dealiased algorithm and rotational algorithms, these sums are slightly negative, due to the fact that wavevectors outside the $k_\text{max}$-radius sphere have not been included.
As expected, the sum of the helicity transfers in the skew-symmetric simulation is significantly different than zero, and in particular it is negative, meaning that helicity is being numerically dissipated.
The nonlinear transfers $T_e(k)$ and $T_h(k)$ are plotted in Fig.~\ref{fig:Tspectra}.
For clarity, only the positive values in the unforced range are shown. 
As expected, for the energy transfers the error is low, the rotational form showing larger values at the largest scales (close to the forcing scale) as well as at the dissipative scales.
Importantly, the skew-symmetric helicity transfer almost superposes with the reference one at larger scales, while it becomes significantly lower than the reference computation at smaller scales. 
The rotational simulation shows helicity transfer values larger than the exact spectrum over almost all the resolved range, and only at the dissipative scales does $T_h(k)$ become smaller than in the de-aliased spectrally truncated case.

As a result, while for the skew-symmetric form the discretization error affects mainly the small scales (causing a numerical helicity dissipation) and leaves the helicity budget basically unaltered at large scales, for the rotational form the nonlinear transfer $T_h(k)$ results to be significantly \textit{distorted} by the discretization error, leading to an increased helicity cascade, although its integral over all scales vanishes.
These observations explain the good behaviour of the skew-symmetric algorithm and the poor behaviour of the rotational algorithm in reproducing the exact helicity spectrum.

\section{Conclusions} \label{sec:conclusions}

The Navier-Stokes equations possess two inviscid quadratic invariants:
the mean kinetic energy $\langle \vv{u} \cdot \vv{u} \rangle/2$ and the mean helicity
$\langle \vv{u} \cdot \vv{\omega} \rangle$. While the former is intimately related to the conservative
structure of the Euler equations, the latter comes instead from the
degeneracies of the Euler operator and, like energy, plays a fundamental role
in turbulent flows.

Upon discretization, the invariant character of the aforementioned
quantities is generally lost due to numerical errors, and particularly
due to the lack of the product rule on a discrete level and the
contamination of truncation and aliasing errors.
Numerical methods capable of preserving energy and helicity might be highly
desirable for reliable computer simulations of turbulent flows, although
the impact and interference of numerical errors have also to be
taken into consideration.

We have characterized the helicity conservation properties of commonly
used discretization algorithms for the Navier-Stokes equations, with particular
reference to the formulation employed to express the convective term.
It is found that, for a colocated layout, spatial conservation of helicity
(or, in the viscous case, stric enforcement of the discrete helicity balance)
is possible through the use of the rotational form. On the other hand, the energy-preserving
and widely used skew-symmetric form fails to preserve helicity.
As far as the time integration is concerned, implicit symplectic methods are required
to satisfy conservation of both quadratic invariants to machine accuracy, while explicit pseudo-symplectic
schemes can ensure higher-order conservation.

Numerical simulations have been carried out in both inviscid and viscous cases.
The inviscid cases have fully confirmed the theoretical results and demonstrated
that the skew-symmetric form dissipates the initial helicity content in few characteristic times.
On the other hand, the rotational form was able to conserve both energy and helicity, 
although its behaviour was found to be strongly dependent on the spatial accuracy of the
schemes (i.e., the truncation error). When used in conjunction with a second-order method,
the rotational form displayed a significant error
with respect to the exact energy and helicity spectra predicted by Kraichnan
at absolute equilibrium.
The inviscid results have also shown that the popular staggered Harlow-Welch scheme does not
preserve helicity, although the definition of this quantity is somewhat
ambiguous on a staggered grid. 
As an additional result, it was found that explicit Runge-Kutta schemes,
which are typically dissipative for energy, can instead be productive for helicity.

In the viscous case, we performed finite-difference forced simulations of helical turbulence at moderate Reynolds number. The simulations have highlighted that, despite the
global conservation properties of the rotational form, the results obtained with
the skew-symmetric form are overall in better agreement with the reference dealiased
spectral simulation. The explanation for this counter-intuitive behaviour 
was obtained by performing a scale-by-scale analysis of the nonlinear energy 
and helicity transfers. It was found that the helicity-conservation error 
of the skew-symmetric form is primarily biased towards the smaller scales.
In particular, a depletion of the helicity transfer is observed, which causes a helicity dissipation 
affecting only the small-scale part of the helicity spectrum.
On the other hand, the nonlinear transfer produced
by the rotational form was found to be significantly distorted, 
leading to helicity spectra in significant disagreement with the reference solution.

\bibliographystyle{spmpsci}      % mathematics and physical sciences
\bibliography{bibliography}   % name your BibTeX data base

\end{document}